\documentclass[12pt,reqno]{article}
\usepackage{amsmath,amsthm,amssymb}

\numberwithin{equation}{section}

\usepackage{tikz-cd}
\usepackage[all]{xy}
\usepackage{mathtext}
\usepackage[T1,T2A]{fontenc}
\usepackage[utf8]{inputenc}
\usepackage[russian,english]{babel}
\usepackage{amsfonts,longtable,amssymb,amsmath,latexsym, dsfont}
\usepackage{wrapfig}
\usepackage{xypic}
\usepackage{mathtools}
\usepackage{hyperref}

\title{Statistical Lie algebras of a constant curvature and locally conformally K\"ahler Lie algebras}
\author{Pavel Osipov\footnote{National Research University Higher School of Economics, Russian Federation }\footnote{Pavel Osipov is 
		partially supported by the HSE University Basic
		Research Program, Simons Foundation, and by the contest “Young Russian
		Mathematics”.} }

\usepackage{amsthm, amsmath, amssymb, amsfonts, amscd, amscd, geometry}
\newsavebox{\ssa}
\renewcommand{\d}{\partial}

\newcommand{\dr}{\frac{\partial}{\partial r}}

\newcommand{\g}{\mathfrak{g}}

\newcommand{\R}{\mathbb R}

\newcommand{\CC}{\mathbb{C}}
\newcommand{\ds}{\frac{\partial}{\partial s}}

\theoremstyle{definition}
\newtheorem{theorem}{Theorem}[section]
\newtheorem{lemma}[theorem]{Lemma}
\newtheorem{proposition}[theorem]{Proposition}

\newtheorem{cor}[theorem]{Corollary}
\newtheorem{defin}[theorem]{Definition}
\newtheorem{example}[theorem]{Example}

\renewcommand{\L}{\mathcal{L}}
\theoremstyle{remark}
\newtheorem{rem}{Remark}[section]
\usepackage{indentfirst}

\begin{document}
	
		\maketitle
		
		\begin{abstract}
			We show that a statistical manifold manifold of a constant non-zero curvature can be realised as a level line of Hessian potential on a Hessian cone. We construct a Sasakian structure on $TM\times\R$ by a statistical manifold manifold of a constant non-zero curvature on $M$. By a statistical Lie algebra of a constant non-zero Lie algebra we construct a l.c.K Lie algebra.

		\end{abstract}
	
		\tableofcontents

	\section{Introduction}
	A {\bfseries flat affine manifold} is a differentiable manifold equipped with a flat, torsion-free connection. Equivalently, it is a manifold equipped with an atlas such that all translation maps between charts are affine transformations (see \cite{FGH} or \cite{shima}). 
	A {\bfseries Hessian manifold} is a flat affine manifold with a Riemannian metric which is locally equivalent to a Hessian of a function. Equivalently, a Hessian manifold is a flat affine manifold $(M,\nabla)$ endowed with Riemnnian metric $g$ such that the tensor $\nabla g$ is totally symmetric. Any Kähler metric can be defined as a complex Hessian of a plurisubharmonic function. Thus, the Hessian geometry is a real analogue of the Kähler one. 
	
	A Kähler structure $(I,g^{\text{r}})$ on $TM$ can be constructed by a Hessian structure $(\nabla,g)$ on $M$ (see \cite{shima}). The correspondence 
	$$
	\text{r}:\{\text{Hessian manifolds}\} \to \{\text{K\"ahler manifolds}\}
	$$
	$$
\ \ \ \	\ (M,\nabla,g)\  \ \to \ \ (TM,I,g^{\text{r}})
	$$
	is called the {\bfseries r-map}. In particular, this map associates special Kähler manifolds to special real manifolds (see \cite {AC}). In this case, r-map describes a correspondence between the scalar geometries for supersymmetric theories in dimension {5 and 4.} See \cite{CMMS} for details on the r-map and supersymmetry.
	
	Hessian manifolds have many different application: in supersymmetry (\cite{CMMS}, \cite{CM}, \cite{AC}), in convex programming
	(\cite{N}, \cite{NN}), in the Monge-Ampère Equation (\cite{F1}, \cite{F2}, \cite{G}), in the WDVV equations (\cite{T}).

	A {\bfseries Riemannian cone} is a Riemannian manifold $(M\times \R^{>0}, s^2g_M+ds^2)$, where $t$ is a coordinate on $\R^{>0}$ and $g_M$ is a Riemannian metric on $M$. Riemannian cones have important applications in supegravity (\cite{ACDM}, 
	\cite{ACM}, \cite{CDM}, \cite{CDMV}). Geometry and holonomy of pseudo-Riemannian cones are studied in \cite{ACGL} and \cite{ACL}.

	A Riemannian manifold $(M,g)$ is called {\bfseries Sasakian} if there exists a complex structure $I$ on the cone $\left(M\times \R^{>0},s^2g_M+ds^2\right)$ such that $\left(M\times \R^{>0},s^2g_M+ds^2,I\right)$ is a K\"ahler manifold. Note that our definitions of Sasakian manifolds are not standard but equivalent. See \cite{5sasaki} or \cite{BG} for standard definitions and \cite{ACHK} or \cite{OV} for equivalence of them.
		
			A {\bfseries radiant manifold} $(C,\nabla, \rho)$ is a flat affine manifold $(C,\nabla)$ endowed with a vector field $\rho$ satisfying. Equivalently, it is a manifold equipped with an atlas such that all translation maps between charts are linear transformations (see e.g. \cite{Go}). A {\bfseries Hessian cone} is a Hessian manifold $(M\times\R^{>0},\nabla,g=s^2g_M+ds^2)$ such that there exists a constant $\lambda 
			\ne 0,\frac{1}{2}$ such that $(M\times \R^{>0},\nabla,\lambda s\ds)$ is a radiant manifold. 
			In the case $\lambda=\frac{1}{2}$, the metric $g$ satisfies 
				$
				\iota_\xi g = 0,
				$
				i.e. can not be positive definite (\cite{G-A}).
			
		A {\bfseries statistical manifold} $\left(M,D,g\right)$ is a manifold $M$ endowed with a torsion-free connection $D$ and a Riemannian metric $g$ such that the tensor $D g$ is totally symmetric. A statistical manifold $\left(M,D,g\right)$ is said to be {\bfseries of a constant curvature} $c$ if the curvature tensor $\Theta_D$ satisfies 
		$$
		\Theta_D (X,Y)Z=c\left(g(Y,Z) X - g(X,Z) Y\right),
		$$ 
		for any $X,Y,Z\in TM$ (see \cite{shima} or \cite{K}). Note that the set of statistical manifolds of constant curvature zero are exactly the set of Hessian manifolds.

		We show that a statistical manifold of a constant curvature can be realised as a level line of a Hessian potential on a Hessian cone. A Sasakian manifold is a level set of K\"ahler cone. In this sense, statistical manifolds of constant curvature are real analogue of Sasakian manifolds.
			\begin{theorem}
				Let $(M, g,\nabla)$ be a statistical manifold of a constant curvature. Then $TM\times\R$ admits a structure of a Sasakian manifold.
			\end{theorem}
	This theorem is closely related to the r-map. Namely, we have a diagram 
	$$
	\begin{CD}
	M\times \R^{>0} @>\text{r}>> T(M\times\R^{>0})\\
	@AA con A @AA con A @.\\
	M @>>> TM\times\R
	\end{CD},
	$$
	where vertical arrows associate Riemannian cones to the corresponding Riemannian manifolds. The theorem implies that the Riemannian manifold $T(M\times\R^{>0})$ with the metric constructed by r-map is actually a cone over $TM\times \R$.

	Then we work with Lie algebras and groups equipped with invariant structures on them. There are different descriptions of an {\bfseries affine structure} on a Lie algebra $\mathfrak{g}$: a torsion-free flat connection on $\mathfrak{g}$, an {\bfseries étale affine representation} $\mathfrak{g}\to \mathfrak{aff}(\R^n)$, where $n$ is dimension of $\mathfrak{g}$, or a structure of a left symmetric algebra on $\mathfrak{g}$, that is, a multiplication on $\mathfrak{g}$ satisfying 
		$$
		XY-YX=[X,Y] \ \ \ \text{and} \ \ \ 
		X(YZ)-(XY)Z=Z(XY)-(ZX)Y
		$$
		for any $X,Y,Z \in \mathfrak{g}$ (see \cite{burde} or \cite{Bu2}).

		An almost complex structure $I$ on a Lie algebra $\g$ is called {\bfseries integrable} if the Nijenhuis tensor of $I$ equals to zero i.e. for any $X,Y\in \g$ 
		$$
		N_I(X,Y)=[X,Y]+I([IX,Y]+[X,IY]-[IX,IY]=0.
		$$
		An almost complex structure on a Lie algebra $\g$ of left invariant fields of a Lie group $G$ sets a left invariant almost complex structure on $G$. It is follows from Newlander–Nirenberg theorem that the almost complex structure on $I$ is integrable if and only if the left invariant almost complex structure on $G$ is integrable.

		Let $(\g,\nabla)$ be a Lie algebra with a flat torsion free connection and $\g_a$ the abelian Lie algebra which coincides with $\g$ as a vector space.  Consider the Lie algebra $\g\ltimes_\nabla \g_{a}$ that is the vector space $\g\oplus \g$ with the commutator
		$$
		[X_1\oplus Y_1, X_2\oplus Y_2] = [X_1,X_2] \oplus\left( \nabla_{X_1} Y_2- \nabla_{X_2} Y_1\right).
		$$
		Then the almost complex structure $I$ on $\g\ltimes_\nabla \g_a$ defined by the rule 
		$$
		I\left(X_1\oplus X_2\right)=-X_2\oplus X_1.
		$$
		is integrable (see \cite{CO} or \cite{BD}). If $\theta$ is an étale affine representation of $G$ then the Lie algebra of left invariant fields on $G\ltimes_\theta \R^n$ equals 
		$\g\ltimes_\nabla \g_a$. 
		
			A {\bfseries Hessian Lie algebra} $\left(\g,\nabla,g\right)$ is a Lie algebra $\g$ endowed with a flat torsion free connection $\nabla$ and symmetric bilinear form $g$ such that $\nabla g$ is totally linear i.e. for any $X,Y,Z\in \g$ we have 
		$$
		g(\nabla_X Y,Z)+g(Y,\nabla_X Z)=g(\nabla_Y X, Z)+g(X, \nabla_Y Z).
		$$
		A {\bfseries Hessian Lie group} $(G,\nabla, g)$ is a Lie group $G$ endowed with a left invariant affine structure $\nabla$ and a left invariant Hessian metric $g$. A Lie groups admits a Hessian structure if and only if the corresponding Lie algebra admits a Hessian structure.

		We adapt r-map to the case of Lie algebras and groups. A {\bfseries Kähler Lie algebra} is a Lie algebra endowed with an integrable almost complex structure $I$ and closed 2-form $\omega$ such that the bilinear form $\omega(\cdot,I\cdot)$ is positive definite.   
				\begin{theorem}\label{T11}
			Let $\left(\g,\nabla,g\right)$ be an Hessian Lie algebra and $\pi:\g\ltimes_\nabla \g_a\to\g$ the projection. Then $\left(\g\ltimes_\nabla \g_a,I,\omega\right)$ is a K\"ahler Lie algebra, where 
			$$
			I(X\oplus Y)=- Y\oplus X \ \ \ \text{and} \ \ \ \omega(X,Y)=\pi^*g(IX,Y)-\pi^*g(X,IY).
			$$
		\end{theorem}

		\begin{cor}\label{T22}
			Let $G$ be an $n$-dimensional simply connected Lie group equipped with a left invariant affine structure $\nabla$ and $\theta$ the linear part of the corresponding affine action of $G$. Then there exists a left invariant Kähler metric on $G\ltimes_\theta \R^n.$
		\end{cor}
	\begin{rem}
			Note that a Kähler structure on a the group $G\ltimes_{\theta^*} \left(\R^n\right)^*$ is constructed by an invariant Hessian structure on $G$ in \cite{medina}. The corresponding complex structure on $G\ltimes_{\theta^*} \left(\R^n\right)^*$ depends on the Hessian metric on $G$. In our case, the complex structure on $G\ltimes_\theta \R^{n}$ depends only on the affine structure on $G$. 
	\end{rem}

		A {\bfseries locally conformally K\"ahler (l.c.K) manifold / Lie algebra} is a manifold / Lie algebra endowed with a complex structure $I$, closed 1-form $\theta$, and 2-form $\omega$ such that 
		$
		d\omega=\theta\wedge \omega
		$ 
		and $\omega(\cdot,I\cdot)$ is positive definite. 
	The closed form $\theta$ on a l.c.K. manifold is locally exact i.e. equals to a differential of a locally defined function $f$. Then the locally defined form $e^{-f}\omega$ is K\"ahler.

			A {\bfseries statistical Lie algebra} $\left(\g,g,D\right)$ is a Lie algebra endowed with a bilinear symmetric positive-definite form $g$ and torsion-free connection and $D$ such that $D g$ is a totally symmetric tensor. A statistical Lie algebra $\left(\g,g,D\right)$ is said to be {\bfseries of a constant curvature} $c$ if the curvature tensor equals 
			$$
			\Theta_D (X,Y)Z=c\left(g(Y,Z) X - g(X,Z) Y\right),
			$$ 
			for any $X,Y,Z\in\g$.

			A {\bfseries statistical Lie group (of a constant curvature $c$)} $\left(G,g,D\right)$ is a Lie group endowed with a left invariant statistical structure (of a constant curvature $c$).

		Obviously, there is one-to-one correspondence between simply connected statistical Lie groups (of a constant curvature $c$) and statistical Lie algebras (of a constant curvature $c$).

    	\begin{theorem}\label{TT}
    	Let $\left(\g,g_\g,D\right)$ be a statistical Lie algebra of a constant nonzero curvature $c$, $\rho$ a generator of the subalgebra $\{0\}\times \R \subset \g\times \R$, $\nabla$ a connection on Lie algebra defined by 
    	$$
    	\nabla_X Y = D_X Y -c g(X,Y)\rho, \ \ \  \nabla_X \rho=\nabla_\rho X=X, \ \ \ \ \nabla_\rho \rho= \rho
    	$$
    	for any $X,Y\in \g\times \{0\}\subset \g\times \R$, $\pi: \left(\g\times \R\right)\ltimes_\nabla\left(\g\times \R\right)_a\to \g$ the projection on the first factor. Denote
    	$
    	\rho_1=\rho\oplus 0, \ \rho_2=0\oplus \rho  \in  \left(\g\times \R\right)\ltimes_\nabla\left(\g\times \R\right)_a.
    	$
    	Consider an almost complex structure $I$ and 2-form $\omega$ on $\left(\g\times \R\right)\ltimes_\nabla\left(\g\times \R\right)_a$ defined by
    		$$
    	I(X\oplus Y)=- Y\oplus X,  \ \ \ \  \omega(X,Y)=\pi^*g(IX,Y)-\pi^*g(X,IY),
    	$$
    	for any $X,Y\in \g\times\R$. 
    	Then for any $t\in \R^{>0}$
    	$$
    	\left(\left(\g\times \R\right)\ltimes_\nabla\left(\g\times \R\right)_a,I,\omega_t=\omega + t\rho^1\wedge\rho^2,-(1+ct)\rho_1^*\right)
    	$$
    	is a l.c.K. Lie algebra. Moreover, if $1+ct=0$ then this algebra is K\"ahler. 
    \end{theorem}
	If $1+ct=0$ then collection $\left(\g,\nabla,g_t=g+t\left(\rho^*\right)^2\right)$ is a Hessian Lie algebra (see \cite{MS} or \cite{shima}). The construction of a K\"ahler Lie algebra from Theorem \ref{TT} arises from applying Theorem \ref{T11} to the Hessian Lie algebra $\left(\g,\nabla,g_t=g+t\left(\rho^*\right)^2\right)$.

    	\begin{cor}\label{1.5}
    	Let $G$ be an $n$-dimensional simply connected statistical Lie group of constant curvature $c$ and $\theta$ the linear part of the corresponding affine representation of $G \times \R^{>0}$. Then there exists an étale affine representation $\mathfrak{g}\to \mathfrak{aff}\left(\R^{n+1}\right)$ and a left invariant l.c.K. structure on  ${\left(G\times\R^{>0}\right)\ltimes_\theta \R^{n+1}}$. Moreover, if $c<0$ then there exists a K\"ahler structure on ${\left(G\times\R^{>0}\right)\ltimes_\theta \R^{n+1}}$.
    \end{cor}

 Further, we provide examples of statistical Lie algebras of a constant curvature and apply the construction of l.c.K Lie algebras to them. First, we consider examples of statistical Lie algebras of a constant curvature called clans. Clans are a curtain class of Lie algebras that can be constructed by a homogeneous cones without full straight lines (see \cite{vinb}). Second, we consider Lie algebras $\mathfrak{so}(2)$ and $\mathfrak{su(2)}$.

	\section{Geometric structures on manifolds}
	
	\subsection{Hessian and Kähler structures}\label{rm}
	\begin{defin}
		A {\bfseries flat affine manifold} is a differentiable manifold equipped with a flat, torsion-free connection. Equivalently, it is a manifold equipped with an atlas such that all translation maps between charts are affine transformations (see \cite{FGH}).  
	\end{defin}
	
	\begin{defin}
		A Riemannian metric $g$ on a flat affine manifold $(M,\nabla)$ is called to be a {\bfseries Hessian metric} if $g$ is locally expressed by a Hessian of a function
		$$
		g=\text{Hess} \ \varphi =\frac{\partial^2\varphi}{\partial x^i \partial x^j} dx^i dx^j,
		$$
		where $x^1,\ldots, x^n$  are flat local coordinates. Equivalently, $g$ is Hessian if and only if the 3-tensor $\nabla g$ is totally symmetric. A {\bfseries Hessian manifold} $(M,\nabla,g))$ is a flat affine manifold $(M,\nabla)$ endowed with a Hessian metric $g$. (see \cite{shima}). 
	\end{defin}
	
	Let $U$ be an open chart on a flat affine manifold $M$, functions $x^1,\ldots, x^n$ be affine coordinates on $U$, and $x^1,\ldots, x^n, y^1, \ldots, y^n$ be the corresponding coordinates on $TU$. Define the complex structure $I$ by $I(\frac{\partial}{\partial x^i})=\frac {\partial} {\partial y^i}$. Corresponding complex coordinates are given by $z^i=x^i+\sqrt {-1}y^i$. The complex structure $I$ does not depend on a choice of flat coordinates on $U$. Thus, in this way, we get a complex structure on the $TM$. 
	
	Let $\pi : TM \to M$ be a natural projection. Consider a Riemannian metric $g$ on $M$ given locally by
	\begin{equation*}\label{1}
	g_{i,j} dx^idx^j.
	\end{equation*}
	
	Define a bilinear form $g^{\text{r}}$ on $TM$ by
	\begin{equation*}\label{2}
	g^{\text{r}}=\pi^* g_{i,j} \left(dx^idx^j+dy^idy^j\right)
	\end{equation*}
	or, equivalently,
	\begin{equation}\label{101}
	g^{\text{r}}(X,Y)=\left(\pi^*g\right)(X,Y)+\left(\pi^*g\right)(IX,IY),
	\end{equation}
	for any $X,Y\in T\left(TM\right)$.
	
	\begin{proposition}[\cite{shima}, \cite{AC}] \label{2.3}\label{23}
		Let $M$ be a flat affine manifold, $g$ and $g^{\text{r}}$ as above. Then the following conditions are equivalent:
		\begin{itemize}
			\item[(i)]$g$ is a Hessian metric.
			
			\item[(ii)] $g^{\text{r}}$ is a Kähler metric. 
		\end{itemize}
		Moreover, if
		$
		g=\text{Hess} \varphi
		$ 
		locally then $g^{\text{r}}$ is equal to a complex Hessian
		$$
		g^{\text{r}}=\text{Hess}_\mathbb{C} (4\pi^*\varphi)=\d\bar{\d}(4\pi^*\varphi).
		$$
	\end{proposition}

	\begin{defin}
		The metric $g^{\text{r}}$ is called a {\bfseries Kähler metric associated to $g$}. The correspondence which associates the Kähler manifold $(TM,g^{\text{r}})$ to a Hessian manifold $(M, g)$ is called the {\bfseries (affine) r-map} (see \cite{AC}).
	\end{defin}
\subsection{Hessian cones and statistical manifolds of a constant nonzero curvature}

\begin{defin}
	A {\bfseries radiant manifold} $(C,\nabla, \rho)$ is a flat affine manifold $(C,\nabla)$ endowed with a vector field $\rho$ satisfying
	\begin{equation}
	\nabla \rho =\text{Id}.
	\end{equation}
	
	Equivalently, it is a manifold equipped with an atlas such that all translation maps between charts are linear transformations (see e.g. \cite{Go}).
\end{defin}

\begin{proposition}[\cite{Go}]\label{Go}\label{32}
	Let $s$ be a coordinate on $\R^{>0}$ and $(M\times \R^{>0},\nabla,s\ds)$ a radiant manifold. Consider a natural action of $\R^{>0}$ on $M\times \R^{>0}$. Then the connection $\nabla$ is $\R^{>0}\text{-invariant}$.
\end{proposition}

\begin{defin}
	A {\bfseries Hessian cone} is a Hessian manifold $(M\times\R^{>0},\nabla,g=s^2g_M+ds^2)$ such that there exists a constant $\lambda\ne 0,\frac{1}{2}$ such that $(M\times \R^{>0},\nabla,\lambda s\ds)$ is a radiant manifold.
\end{defin}

\begin{proposition}[\cite{G-A}]\label{GA}
	Let $(C,\nabla, \rho)$ be a radiant manifold and $g$ a Hessian metric on $M$ with respect to $\nabla$. Then 
	$$
	\L_\rho  g =g +\nabla(\iota_\rho g).
	$$
\end{proposition}
\begin{proposition}\label{29}
	Let $(M\times\R^{>0},\nabla,g=s^2g_M+ds^2)$ is a Hessian cone. Then we have
	\begin{equation}\label{eh}
	g=\text{Hess} \left(\frac{\lambda s^2}{4\lambda-2}\right) .
		\end{equation}
\end{proposition}
\begin{proof}
	Let $\rho=\lambda s\ds$ be the radiant vector field. Then
	$$
	\iota_\rho g = \iota_{\lambda s\ds} \left(s^2g_M+ds^2\right)= \lambda  sds=d\left(\frac{\lambda s^2}{2}\right).
	$$
	Hence,
	$$
	\text{Hess}\left(\frac{\lambda  s^2}{2}\right)=\nabla d\left(\frac{\lambda  s^2}{2}\right)=\nabla\iota_{\rho} g.
	$$
	Combining this with Proposition \ref{GA}, we get 
	$$
	\L_\rho g - g = \text{Hess}\left(\frac{\lambda s^2}{2}\right).
	$$
	We have $\L_\rho g = \L_{\lambda r\dr} \left(s^2g_M+ds^2\right)=2\lambda g$.
	Thus,
	$$
	g=\text{Hess} \left(\frac{\lambda s^2}{4\lambda -2}\right)  .
	$$
\end{proof}

\

\begin{defin}
	A {\bfseries statistical manifold} $\left(M,D,g\right)$ is a manifold $M$ endowed with a torsion-free connection $D$ and a Riemannian metric $g$ such that the tensor $D g$ is totally symmetric. A statistical manifold $\left(M,D,g\right)$ is said to be {\bfseries of a constant curvature} $c$ if the curvature tensor $\Theta_D$ satisfies 
	$$
	\Theta_D (X,Y)Z=c\left(g(Y,Z) X - g(X,Z) Y\right),
	$$ 
	for any $X,Y,Z\in TM$.
\end{defin}


\begin{lemma}\label{l}
	Let $(M\times\R^{>0},\nabla,g=s^2g_M+ds^2)$ be a Hessian cone and $\rho=\lambda s\ds$ the radiant vector field. Then the following condition are satisfies:
	\begin{itemize}
		\item [(i)] $
		\left(\nabla_X g\right)\left(Y, \rho\right)=\left(\nabla_\rho g\right)\left(X, Y\right)=\left(2\lambda -2\right)g(X,Y).
		$
		\item [(ii)]  $g\left(\nabla_X Y, \rho\right)=\left(1-2\lambda \right)g(X,Y)$
	\end{itemize} 
\end{lemma}
\begin{proof}
	(i) Since $g$ is a Hessian metric, the tensor $\nabla g$ is totally symmetric. Therefore, 
	\begin{equation}\label{35}
	\left(\nabla_X g\right)\left(Y, \rho\right)=\left(\nabla_\rho g\right)\left(X, Y\right)=\L_\rho\left(g(X,Y)\right)-g(\nabla_\rho X,Y)-g(X,\nabla_\rho Y).
	\end{equation}
	Since $X,Y\in TM$, the value $g_M(X,Y)$ is constant along $\rho$. Hence, 
	$$
	\L_\rho \left(g(X,Y)\right)= \L_\rho \left(s^2g_M(X,Y)\right)=2 \lambda \left(s^2g_M(X,Y)\right)=2\lambda  \left(g(X,Y)\right).
	$$
	Moreover, $\rho$ commutes with $X$ and $Y$. Hence, $\nabla_\rho X=\nabla_X\rho=X$ and $\nabla_\rho Y=\nabla_Y\rho=Y$. Combining this with \eqref{35}, we get $
	\left(\nabla_X g\right)\left(Y, \rho\right)=\left(\nabla_\rho g\right)\left(X, Y\right)=\left(2\lambda -2\right)g(X,Y).
	$

	(ii) We have 
	$$
	g\left(\nabla_X Y, \rho\right)=-\left(\nabla_X g\right)\left(Y, \rho\right)+\L_X \left(g\left(Y, \rho\right)\right)-g\left( Y, \nabla_X \rho\right)= -\left(\nabla_X g\right)\left(Y, \rho\right)-g(Y,X).
	$$
	Combining this with the item (i), we obtain
	$
	g\left(\nabla_X Y, \rho\right)=\left(1-2\lambda \right)g(X,Y).
	$
\end{proof}

\begin{theorem}\label{T}
	
	Let $(M\times\R^{>0},\nabla,g=s^2g_M+ds^2)$ be a Hessian cone, $\rho=\lambda s\ds$ the radiant vector field, and $c=\frac{2\lambda -1}{\lambda ^2}$. Let us identify $M$ with the submanifold $M\times 1\subset M\times\R^{>0}$.  Then for any $X,Y\in TM$, we have 
	\begin{equation}\label{AA}
	\nabla_X Y = D_X Y -c g_M(X,Y)\rho,
	\end{equation}
	where $D$ is a torsion-free connection on $M$. Moreover, $\left(M,g_M,D\right)$ is a statistical manifold of curvature $c$.
	
	Conversely, if $\left(M,g_M,D\right)$ is a statistical manifold of a non-zero constant curvature $c\le 1$, $\lambda $ a solving of the equation $\frac{2\lambda -1}{\lambda ^2}=c$, and $\rho=\lambda  s\ds$ a field of on $M\times \R^{>0}$ then $(M\times\R^{>0},\nabla,g=s^2g_M+ds^2)$ is a Hessian cone, where the connection $\nabla$ is defined by equation \eqref{AA} and
	\begin{equation}\label{AAAA}
	\nabla_X \rho=\nabla_\rho X=X, \ \ \ \ \nabla_\rho \rho= \rho,
	\end{equation}
	for any $X\in TM$.

\begin{proof}
	Since $\rho$ is orhtogonal to $M$ we have 
	$$
	\nabla_X Y= D_X Y+ \frac{g(\nabla_X Y, \rho)}{g(\rho,\rho)} \rho=D_XY+\frac{g(\nabla_X Y, \rho)}{\lambda  ^2s^2}\rho.
	$$
	where $D_X Y\in TM$. Combining this with item (ii) of Lemma \ref{l}, we get that 
	$$
	\nabla_X Y= D_X Y+ \frac{g(\nabla_X Y,\rho)}{\lambda  ^2s^{2}} \rho=D_X Y +\frac{1-2\lambda }{\lambda ^2s^2}g(X,Y)\rho.
	$$
	Combining this with $g(X,Y)=s^2g_M(X,Y)$ we get \eqref{AA}. The term $g(*,*) \rho$ is 2-1 tensor and $\nabla$ is a connection. Therefore, $D$ is a connection.
	
	The connection $\nabla$ is flat hence for any $X,Y,Z\in TM$ we have
	$$
	\Theta_\nabla(X,Y) Z=\nabla_X\nabla_YZ-\nabla_Y\nabla_XZ-\nabla_{[X,Y]}Z=0.
	$$
	Combining this with \eqref{AA} and $\nabla\rho=\text{Id}$ we get
	\begin{multline*}
	\Theta_\nabla(X,Y) Z=\nabla_X\nabla_YZ-\nabla_Y\nabla_XZ-\nabla_{[X,Y]}Z=\Theta_D(X,Y) Z-c \left(g_M(Y,Z) X- g_M(X,Z) Y\right)+\\
	-c \left(\L_Y\left(g_M(X,Z)\right)-\L_X\left(g_M(Y,Z)\right)+g_M(Y,D_X Z)-g_M(X,D_Y Z)+g_M([X,Y],Z)\right) \rho.
	\end{multline*}
	Since $X,Y,Z\in TM$ and $\Theta_\nabla=0$, we get that
	$$
	\Theta_D(X,Y) Z= c \left(g_M(Y,Z) X -g_M(X,Z) Y\right)
	$$
	and
	$$
	\L_X\left(g_M(Y,Z)\right)- \L_Y\left(g_M(X,Z)\right)+g_M(X,D_Y Z)-g_M(Y,D_X Z)-g_M([X,Y],Z)=0.
	$$
	Combining the last equation with the identity $[X,Y]=D_X Y-D_Y X$ and the formula of covariant derivative of a metric we get that
	$$
	\left(D_X g\right)(Y,Z)-\left(D_Y g\right)(X,Z)=0.
	$$
	Thus, the tensor $Dg$ is totally symmetric. We proved the first part of the theorem.
	
	Now, let $\left(M,g_M,D\right)$ be a statistical manifold of curvature $c$. Then $\nabla$ is flat by the same calculation as above (in the opposite side). Thus, it is enough to check that the metric $g$ is Hessian. 
	
	For any $X,Y,Z \in TM$ we have 
	\begin{multline*}\label{AAA}
	\left(\nabla_X g\right)(Y,Z)=\L_X\left(g(Y,Z)\right)-g\left(\nabla_X Y,Z\right)-g\left(Y,\nabla_X Z\right)= \\
	=\L_X\left(g(Y,Z)\right)-g\left(D_X Y,Z\right)-g\left(Y,D_X Z\right)=\left(D_X g\right)(Y,Z)
	\end{multline*}
	Combining this with $\left(D_X g\right)(Y,Z)-\left(D_Y g\right)(X,Z)=0$, we get 
	$$
	\left(\nabla_X g\right)(Y,Z)=\left(\nabla_Y g\right)(X,Z).
	$$
	We have,
	$$
	\left(\nabla_\rho g\right) (X, Y)=\L_\rho\left(g(X,Y)\right)-g\left(\nabla_\rho X, Y\right)-g\left(X,\nabla_\rho Y\right)=(2a-2)g(X,Y)
	$$
	and
	\begin{multline*}
	\left(\nabla_X g\right)(\rho,Y)=\L_X\left(g(\rho,Y)\right)-g\left(\nabla_X \rho,Y\right)-g\left(\rho,\nabla_X Y\right)=\\ =g(X,Y)-\frac{1-2a}{a^2}g(\rho,\rho)g_M(X,Y)=(2a-2)g(X,Y).
	\end{multline*}
	Thus,
	$$
	\left(\nabla_\rho g\right) (X, Y)=	\left(\nabla_X g\right)(\rho,Y).
	$$
	Finally, 
	$$
	\left(\nabla_\rho g\right)(\rho,X)=\L_\rho\left(g(\rho,X)\right)-g\left(\nabla_\rho \rho,X\right)-g\left(\rho,\nabla_\rho X\right)=0
	$$
	and
	$$
	\left(\nabla_X g\right)(\rho,\rho)=\L_X\left(g(\rho,\rho)\right)-g\left(\nabla_X \rho,\rho\right)-g\left(\rho,\nabla_X \rho\right)=0
	$$
	We checked that the tensor $\nabla g$ is totally symmetric. Therefore, the metric $g$ is Hessian.  
\end{proof}
	\begin{cor}
		Any statistical manifold of a non-zero constant curvature can be realised as a level line of Hessian potential on a Hessian cone.
	\end{cor}
\begin{proof}
	According to Theorem \ref{T}, a statistical manifold of a nonzero constant curvature can be realised as a level of the function $s$ on a Hessian cone $(M\times\R^{>0},\nabla,g=s^2g_M+ds^2)$. It follows from Proposition \ref{29}, that level lines of $s$ coincide with level lines of a Hessian potential.
\end{proof}

\end{theorem}

	\subsection{Statistical manifolds of a nonzero constant curvature and Sasakian manifolds}

	\begin{defin}
		A {\bfseries Sasakian manifold} is a Riemannian manifold $(M,g_M)$ such that the cone metric $s^2 g_M+ ds^2$ on $M\times \R^{>0}$ is Kähler with respect to a dilation invariant complex structure.
	\end{defin}
	
	\begin{proposition}[\cite{OV}]\label{42}
		Let $(M\times \R^{>0}, g, I)$ be a Kähler manifold. For any $q\in \R^{>0}$ consider a map 
		$
		\mu_q : M\times \R^{>0} \to M\times \R^{>0}
		$
		defined by
		$
		\mu_q (m\times s)=m\times qs.
		$ 
		If $\mu_q^* g=q^2 g$ then $g=s^2g_M+ds^2$ and $(M,g_M)$ is a Sasakian manifold. 
	\end{proposition}

	There exists a decomposition	
	$$
	T(M\times\R^{>0})=TM \times T\R^{>0}=TM\times \R \times \R^{>0}.
	$$
	If $M\times\R^{>0}$ possess a Hessian structure then, according to Proposition \ref{23}, $T(M\times\R^{>0})$ admits a Kähler structure.
	
	\begin{proposition}\label{43}
		Let $(M\times \R^{>0}, \nabla, g)$ be a Hessian cone and $g^{\text{r}}$ the constructed by the r-map metric on $T(M\times \R^{>0})$ . Consider $T(M\times\R^{>0})=TM\times \R \times \R^{>0}$ as a cone over $TM\times \R$. Then for any $q\in\R^{>0}$ we have $\eta_q^* g=q^2 g$, where the map 
		$
		\eta_q : TM\times\R\times \R^{>0} \to TM\times \R\times \R^{>0}
		$
		is defined by
		$
		\mu_q (m\times s\times t)=m\times s\times  qt.
		$  
	\end{proposition}
	
	\begin{proof}
		We have the commutative diagram 
		$$
		\begin{CD}
		T(M\times\R^{>0}) @>\eta_q>> T(M\times\R^{>0})\\
		@VV\pi V @VV\pi V  @.\\
		M\times\R^{>0} @>\mu_q>> M\times\R^{>0}
		\end{CD},
		$$
		where $\eta_q$ and $\mu_q$ are multiplications of the coordinate on $\R^{>0}$ by $q$. By the definition of $g^r$ we have 
		\begin{equation}\label{e41}
		g^{\text{r}}(X,Y)=\pi^*g(X,Y)+\pi^*g(IX,IY).
		\end{equation}	
		Since the diagram is commutative, it follows that \begin{equation}\label{e42}
		\eta_q^*\pi^*=\pi^*\mu_q^*.
		\end{equation} 
		Moreover, $g$ is a cone metric. Hence,
		\begin{equation}\label{e43}
		\mu_q^* g= q^2 g.
		\end{equation}
		It follows from \eqref{e41}, \eqref{e42}, and \eqref{e43} that 
		$$
		\eta_q^* g^{\text{r}}(X,Y)=q^2 g^{\text{r}}(X,Y).
		$$

	\end{proof}

	\begin{theorem}\label{T1}
		Let $(M, g_M,D)$ be a statistical manifold of a non-zero constant curvature. Then $TM\times\R$ admits a structure of a Sasakian manifold.
	\end{theorem}
	\begin{proof}
	If the curvature of $(M, g_M,D)$ is $c>1$ then $(M,\frac{1}{c}{g_M}, D)$ is a statistical manifolds of curvature 1. Thus, we can assume that the curvature of $(M, g_M,D)$ does not exceed 1. By Theorem \ref{T}, there exists a Hessian cone $(M\times\R^{>0},\nabla,g=r^2g_M+dr^2)$. Then the r-map defines a K\"ahler structure $(g^{\text{r}},I)$ on $T(M\times\R^{>0})=TM\times\R\times\R^{>0}$. By Proposition \ref{32}, the connection $\nabla$ is $\R^{>0}$ invariant. Therefore, the constructed by $\nabla$ complex structure $I$ is $\R^{>0}$-invariant.  Let us identify $TM\times\R$ with $TM\times\R\times 1\subset TM\times\R\times\R^{>0}$. Combining propositions \ref{42} and \ref{43}, we get that $\left(TM\times \R,g^\text{r}|_{TM\times\R}\right)$ is a Sasakian manifold. 
	\end{proof}

	\section{Geometric structures on Lie groups and algebras}
	\subsection{Flat torsion free connections and complex structure}
		
		The group of affine transformations $\text{Aff} (\mathbb{R}^n)$ is given by the matriсes of the form
		$$
		\begin{pmatrix}
		A & a\\
		0 & 1
		\end{pmatrix}
		\in \text{GL}(\R^{n+1}),
		$$
		where $A\in \text{GL}(\R^n)$, and $a\in \R^n$ is a column vector. The corresponding Lie algebra $\mathfrak{aff}(\R^n)$ is given by matrices of the form
		$$
		\begin{pmatrix}
		A & a\\
		0 & 0
		\end{pmatrix}
		\in \mathfrak{gl}(\R^{n+1}).
		$$
		The commutator of $\mathfrak{aff}(\R^n)$ is equal to
		$$
		\left[
		\begin{pmatrix}
		A & a\\
		0 & 0
		\end{pmatrix},
		\begin{pmatrix}
		B & b\\
		0 & 0
		\end{pmatrix}
		\right]
		=
		\begin{pmatrix}
		[A,B] & A(b)-B(a)\\
		0 & 0
		\end{pmatrix}.
		$$
		Algebra $\mathfrak{aff} (\R^n)$ is the semidirect product $\mathfrak{gl}(\R^n)\ltimes \R^n$, where the commutator is given by
		$$
		[A\ltimes a,B \ltimes b]=[A,B]\ltimes (Ab-Ba).
		$$
		Group $\text{Aff}(\R^n)$ is the semidirect product $\text{GL}(\R^n)\ltimes \R^n$, where multiplication is given by 
		$$
		(A\ltimes a) (B \ltimes B)=AB\ltimes (a+Ab).
		$$

		\begin{defin}
			An affine representation $G\to \text{Aff}(\R^n)$ is called {\bfseries étale} if there exists a point $x \in \R^n$ such that the orbit of $x$ is open and the stabilizer of $x$ is discrete. A representation of a Lie algebra $\g \to  \mathfrak{aff} (\R^n)$ is called {\bfseries étale} if the corresponding representation of the simply connected Lie algebra is étale.
		\end{defin}
		
		\begin{theorem}[\cite{burde} or \cite{Bu2}]\label{b}
			Let $G$ be an $n$-dimensional Lie group and $\g$ the corresponding Lie algebra. Choose an identification $i: \g \to \R^n$. Then there is a 1-1 correspondence between left invariant  flat torsion-free connections and étale affine representations. Moreover, if ${\theta : G\to \text{GL}(\R^n)}$ is the linear part of the étale affine representation corresponding to a left invariant  flat torsion-free connection $\nabla$ then the differential $\eta:\g \to\mathfrak{gl}\left(\R^n\right)$ is defined by $X\to i \circ \nabla_X\circ i^{-1}$.

		\end{theorem}

	
	Let $\nabla$ be a flat connection on a Lie algebra $\g$. Denote by $\g_a$ the vector space $\g$ with the structure of an abelian Lie algebra. Look on $\nabla$ as a Lie algebras representation $\nabla: \g \to  \mathfrak{gl}(\R^n) $, $X\to \nabla_X$. Consider the semidirect product $\g\ltimes_\nabla \g_a$ that is vector space $\g\oplus\g$ with the Lie bracket defined by
	\begin{equation}\label{e}
	[X_1\oplus Y_1, X_2\oplus Y_2] = [X_1,X_2] \oplus\left( \nabla_{X_1} Y_2- \nabla_{X_2} Y_1\right).
	\end{equation}
	The flatness of $\nabla$ is equivalent to the Jacobi identity on  $\g\ltimes_\nabla \g_a$. Define an almost complex structure on $\g\ltimes_\nabla \g_a$ by the rule 
	$$
	I\left(X_1\oplus X_2\right)=-X_2\oplus X_1.
	$$
	\begin{theorem}[\cite{CO} or \cite{BD}]\label{53}
		Let $\nabla$ be a flat connection on a Lie algebra $\g$, and $I$ be an almost complex structure on $\g\ltimes_\nabla \g_a$. Then the almost complex structure $I$ is integrable if and only if $\nabla$ is torsion free.
	\end{theorem}
	\begin{proposition}\label{54}
		Let $\nabla$ be a flat torsion free connection on a Lie algebra $\g$ and $\theta:G\to\text{GL}(\R^n)$ the linear part corresponding étale affine representation. Then the Lie algebra of left invariant fields on $G\ltimes_\theta \R^n$ equals $\g\ltimes_\nabla \g_a$.
	\end{proposition}

	\begin{proof}
		Choosing an identification $i: \g \to \R^n$. According to Theorem \ref{b}, the differential of $\theta$ is defined by $\eta (X)=i \circ \nabla_X\circ i^{-1}$. Then the corresponding lie algebra $g\ltimes_\eta \R^n$ is isomorphic to $\g\ltimes_\nabla\g_a$. 
	\end{proof}	

	\begin{cor}\label{3.5}
		Let $G$ be a simply connected Lie group equipped with a left invariant affine structure, $\mathfrak{g}$ the corresponding Lie algebra, and $\theta$ the linear part of the corresponding affine action of $G$. Then there exists a left invariant integrable complex structure on the group 
		$
		G\ltimes_\theta \R^n.
		$ 
	\end{cor}

	\begin{proof}
		The existence of a left invariant integrable complex structure follows from Theorem \ref{53} and Proposition $\ref{54}$. 
	\end{proof}

	\subsection{Hessian and K\"ahler structures}
	\begin{defin}
		{\bfseries A Hessian Lie group} $(G,\nabla, g)$ is a Lie group $G$ endowed with a left invariant affine structure $\nabla$ and a left invariant Hessian metric $g$. 
	\end{defin}
\begin{defin}
	A {\bfseries Hessian Lie algebra} $\left(\g,\nabla,g\right)$ is a Lie algebra $\g$ endowed with a flat torsion free connection $\nabla$ and symmetric bilinear form $g$ such that $\nabla g$ is totally linear i.e. for any $X,Y,Z\in \g$ we have 
	$$
	g(\nabla_X Y,Z)+g(Y,\nabla_X Z)=g(\nabla_Y X, Z)+g(X, \nabla_Y Z).
	$$
\end{defin}

A Lie groups admits a left invariant Hessian structure if and only if the corresponding Lie algebra admits a Hessian structure.


\begin{proof}[Proof of Theorem \ref{T11}]
	According to Theorem \ref{53} the almost complex structure $I$ is integrable. The bilinear form 	
	$$
	g^\text{r}(X,Y)=\omega(X,IY)=\pi^*g(X,Y)+\pi^*g(IX,IY)
	$$ 
	is positive definite. Hence, it is enough to check that the form $\omega$ is closed.
	For any $X_1,X_2\in \g\oplus 0$ we have
	\begin{equation}\label{ee}
	\omega(X_1\oplus0,X_2\oplus0)=\omega(0\oplus X_1,0\oplus X_2)=0,
	\end{equation}
	Combining this with the formula of exterior derivative
	\begin{equation}\label{3.2}
	d\omega(V_1,V_2,V_3)=-\omega([V_1,V_2],V_3)+\omega([V_1,V_3],V_2)-\omega([V_2,V_3],V_1)
	\end{equation}
	and the definition of the Lie bracket \eqref{e}
	$$
	d\omega(X\oplus 0,Y \oplus 0, Z\oplus 0)=0, \ \ \   d\omega(X\oplus 0, 0 \oplus Y, 0\oplus Z)=0,  \ \ \ d\omega(0\oplus X, 0 \oplus Y, 0\oplus Z)=0.
	$$
    and
	\begin{multline*}
	d\omega(X\oplus 0,Y \oplus 0, 0\oplus Z)=-g([X,Y],Z)-g(Y,\nabla_X Z)+g(X,\nabla_Y Z)= \\
	=-g(\nabla_X Y,Z)-g(Y,\nabla_X Z)+g(\nabla_Y X,Z)+g(X,\nabla_Y Z)=-\left(\nabla_X g\right)(Y,Z)+\left(\nabla_Y g\right)(X,Z)=0.
	\end{multline*}

\end{proof}

Corollary \ref{T22} follows from Proposition \ref{54} and Theorem \ref{T11}. 

\subsection{Statistical structures of a nonzero constant curvature and l.c.K structures}

\begin{defin}
	A {\bfseries statistical Lie algebra} $\left(\g,g,D\right)$ is a Lie algebra endowed with a bilinear symmetric positive-definite form $g$ and torsion-free connection and $D$ such that $D g$ is a totally symmetric tensor. A statistical Lie algebra $\left(\g,g,D\right)$ is said to be {\bfseries of a constant curvature} $c$ if the curvature tensor equals 
	$$
	\Theta_D (X,Y)Z=c\left(g(Y,Z) X - g(X,Z) Y\right),
	$$ 
	for any $X,Y,Z\in\g$.
\end{defin}
\begin{defin}
	A {\bfseries statistical Lie group (of a constant curvature $c$)} $\left(G,g,D\right)$ is a Lie group endowed with a left invariant statistical structure (of a constant curvature $c$). 
\end{defin}

	Obviously, there is one-to-one correspondence between simply connected statistical Lie groups (of a constant curvature $c$) and statistical Lie algebras (of a constant curvature $c$).
	
		\begin{defin}
		A {\bfseries locally conformally K\"ahler (l.c.K) manifold (Lie algebra)} is a manifold (Lie algebra) endowed with an integrable almost complex structure $I$, closed 1-form $\theta$, and 2-form $\omega$ such that
		$$
		d\omega=\theta \wedge \omega.
		$$ 
	\end{defin}

	\begin{proof}[Proof of Theorem \ref{TT}]
	
	The bilinear form $g_t^\text{r}=\omega_t(\cdot,I\cdot)$ equals $g^\text{r}+t\left(\left(\rho^1\right)^2+\left(\rho^1\right)^2\right)$ where $g^\text{r}$ is defined by
	$$
	g(X,Y)=\pi^*g(X,Y)+\pi^*g(IX,IY).
	$$
	For $t>0$ the form $g_t^\text{r}$ is positive definite.
	
	The same computations as in the proof of Theorem \ref{T11} show that for any $X,Y,Z\in \g$ we have
	\begin{equation}
	\omega_t(X\oplus 0, Y\oplus 0, Z\oplus 0)=\omega_t(X\oplus 0, Y\oplus 0, 0\oplus Z)=\omega_t(X\oplus 0, 0\oplus Y, 0\oplus Z)=\omega_t(0\oplus X, 0\oplus Y, 0\oplus Z)=0.
	\end{equation}
	
	Combining the formula of the exterior derivative \eqref{3.2} with \eqref{ee}  and the definition of the Lie bracket \eqref{e} we get that
	for any $X,Y\in \g$ 
	$$
	d\omega_t(0\oplus X,0\oplus Y, 0\oplus \rho)=0, \ \ \ d\omega_t(X\oplus 0,0\oplus Y, 0\oplus \rho)=0,
	$$
	and
	$$
	d\omega_t(X\oplus 0,Y\oplus 0, 0\oplus \rho)=-\omega_t([X,Y]\oplus 0, 0\oplus \rho)+\omega_t(0\oplus X,Y\oplus 0) -\omega_t(0\oplus Y,X\oplus 0)=0.
	$$
	Since $[X,Y]\in \g$, we have $\omega([X,Y]\oplus 0, 0\oplus \rho)=0$. It follows from the definition of $\omega_t$ that
	$$
	\omega_t(0\oplus X, Y\oplus 0)=g(X,Y)=\omega_t(0\oplus Y,X\oplus 0).
	$$
	Therefore, 
	$$
	d\omega_t(X\oplus 0,Y\oplus 0, 0\oplus \rho)=0.
	$$
	Also, using \eqref{ee}, \eqref{3.2}, and \eqref{e} we get that for any $X,Y\in \g$ 
	$$
	d\omega_t(X\oplus 0,Y\oplus 0, \rho\oplus 0)=0, \ \ \ d\omega_t(0\oplus X,0\oplus Y, \rho\oplus 0)=0
	$$
	and
	$$
	g(X\oplus 0,0\oplus Y, \rho\oplus 0)=-\omega_t(0\oplus \nabla_X Y,\rho\oplus 0)-\omega_t(X\oplus 0, 0\oplus Y).
	$$
	We have 
	$$
	\omega_t(0\oplus \nabla_X Y,\rho\oplus 0)=\omega_t (0\oplus D_X Y-cg(X,Y), \rho\oplus 0)=\omega_t(0\oplus -cg(X,Y)\rho,\rho\oplus 0).
	$$
	Combining this with the identity $\omega_t=\omega+t\rho^1\wedge\rho^2$ and the definition of $\omega_t$ we get 
	$$
	\omega_t(0\oplus \nabla_X Y,\rho\oplus 0)=tcg(X,Y)=-tc\omega_t(X\oplus 0,0\oplus Y).
	$$
	Thus, 
	$$
	\omega(0\oplus \nabla_X Y,\rho\oplus 0)=-(1+tc)g(X\oplus0,0\oplus Y).
	$$
	As above, using \eqref{ee}, \eqref{3.2}, and \eqref{e} we get that for any $X\in \g$, we have 
	$$
	\omega_t(X\oplus 0, \rho\oplus 0,0\oplus \rho)=0, \ \ \ \omega_t(0\oplus X, \rho\oplus 0,0\oplus \rho)=0.
	$$
	We checked that
	$$
	d\omega_t=-(1+tc)\rho^1\wedge\omega_t.
	$$
	If the $1+tc=0$. Then the form $\omega_t$ is K\"ahler. 
	
	
	\end{proof}
Corollary \ref{1.5} follows from Proposition \ref{54} and Theorem \ref{TT}.

	\subsection{Examples} 
	\subsubsection{Convex regular cones and clans}\label{7}
				
				\begin{defin}
					A subset $V\subset \R^n$ is called {\bfseries regular} if $V$ does not contain any straight full line.
				\end{defin}

	Let $V\subset \R^n$ be a convex regular domain. We denote the maximal subgroup of $\text{GL}(\R^n)$ preserving $V$ by $\text{Aut}(V)$. Note that if $V$ is a regular convex cone then $$\text{Aut}(V)=(\text{Aut}(V)\cap \text{SL}(\R^n))\times \R^{>0}.$$ 

	The following theorem summarized known results.
	
	\begin{theorem}[\cite{vinb}, \cite{shima}]\label{52}
		Let $V$ be a homogeneous convex cone. Then there exists a function $\varphi$ a subgroup $T\subset\text{Aut}(V)$ satisfying the following conditions.
		\begin{itemize}
			\item[(i)] $T$ acts on $V$ simply transitively.
			\item[(ii)] The bilinear form $g_{con} =\text{Hess}\  \left(\ln \varphi\right)$ is a $T$-invariant Hessian metric.
			\item [(iii)] The group $T_{\text{SL}}=T\cap \text{SL}(\R^n)$ preserves the hypersurface $M=\{\varphi=1\}$ and acts simply transitively on it. 
			\item[(iv)] Let $\pi: TV|_M\to TM$ be the projection along the radiant vector field and $\nabla$ be the standard connection on $\R^n$. Then $\left(D,g_{con}\right)$ be a statistical structure of a constant curvature on $M$, where $D$ is defined by $D_X Y= \pi\left(\nabla_X Y\right)$. 
		\end{itemize}
	\end{theorem}
	
	\begin{defin}
		The hypersurface $M$ from Theorem \ref{52} is called the {\bfseries characteristic hypersurface} of a cone. 
	\end{defin}

		\begin{defin}
		Let $T_{\text{SL}}$ be the Lie groups from Theorem \ref{52}. The corresponding Lie algebra $\mathfrak{t}$ is called a {\bfseries clan}. There exist a purely algebraic definition of clans (see \cite{vinb}).
	\end{defin}

			By Theorem \ref{52}, a characteristic hypersurface $M$ admits $T_{\text{SL}}$-invariant statistical structure of a constant curvature. Thus, any clan admits a statistical structure of a constant curvature.


	 

						\begin{example}
							Let $V$ be the vector space of all real symmetric matrices of rank $n$ and $\Omega$ the set of all positive definite symmetric matrices in $V$. Then $\Omega$ is a regular convex cone and the group of upper triangular matrix $\text{T}(\R^n)$ acts simply transitively on $\Omega$ by $s(x)= s x s^T$, where $x \in \Omega$ and $s \in \text{T}(n)$.The characteristic function is equal to
							$$
							\varphi(x)=(\det x)^{-\frac{n+1}{2}}\varphi (e),
							$$ 
							where $e$ is the unit matrix (see \cite{shima}). The corresponding clan $\mathfrak t$ is the algebra of upper triangular traceless matrices.
							
							Consider the case $n=3$. Here, $\mathfrak t$ is generated by elements
							$$
							u=\begin{pmatrix}
							1 & 0 \\
							0 & -1 \\
							\end{pmatrix}, \ \ \
							v= \begin{pmatrix}
							0 & 1 \\
							0 & 0 \\
							\end{pmatrix}.
							$$
							and the relative $[u,v]=2v$. The statistical structure $\left(D,g\right)$ of costant curvature $-c<0$ is defined by
							$$
							D_u u= D_u v=0, \ \ \ D_v u= -2v, \ \ \ D_v v= u, \ \ \ g_\g=\frac{4 \left(u^*\right)^2+2\left(v^*\right)^2}{c}
							$$
							(see \cite{shima}). 
							
							Then the corresponding l.c.K Lie algebra admits generators $u_1,v_1,\rho_1,u_2,v_2,\rho_2$ and relations
							$$
							[u_1,u_2]=-{4}\rho_2, \ \ \ \ \ \  [v_1,v_2] = -{2} \rho, \ \ \ \ \ \  [\rho_1,\rho_2]=\rho_2, 
							$$
							$$
							[u_1,v_1]=2v_1, \ \ \ [v_1, u_2]=-2v_2, \ \ \ [v_1,v_2]=u_2, \ \ \ [\rho_1,u_2]=u_2, \ \ \ [\rho_1,v_2]=v_2. 
							$$
							The complex structure is defined by
							$$
							I(u_1)=u_2, \ \ \ I(v_1)=v_2, \ \ \  I(\rho_1)=\rho_2. 
							$$
							For any $c,t\in\R^{>0}$ we have the following l.c.K. form  
							$$
							\omega_{c,t} =\frac{4}{c}u^1\wedge u^2+\frac{2}{c}v^1\wedge v^2	+t\rho^1\wedge\rho^2 .
							$$
							In particular, the form 
							$
							\omega_{1,1} ={4}u^1\wedge u^2+2v^1\wedge v^2	+\rho^1\wedge\rho^2 
							$ 
							is K\"ahler.
						\end{example}

	\subsubsection{$\mathfrak{so}(2)$ and $\mathfrak{su}(2)$}\label{8}
	
		\begin{example}
			Consider the group $\R$ as the universal covering of $\text{U}(1)=\text{SO(2)}$. The identification $\text{SO}(2)\times{\R^{>0}} \simeq \R^2 \backslash \{0\}$ sets a Hessian structure $\left(\nabla,g\right)$ on $\R\times{\R^{>0}}$ such that $\lambda_q^* g=q^2g$, where $\lambda_q :\R\times{\R^{>0}}\to \R\times{\R^{>0}}$, $\lambda_q\left(x\times s\right)=x\times qs$. By the same way as in the proof of Theorem \ref{1.5}, we can define a l.c.K structure on the group of homothetic motions of plane $\text{H}(2)=\left(\R\times\R^{>0}\right)\ltimes \R^{2}$. The group $\left(\text{SO}(2)\times{\R^{>0}}\right)\ltimes \R^2$ is equals to the group of the group of homothetic motions of plane $\text{H}(\R^2)$. Thus, we get a l.c.K structure on the univarsal covering $\widetilde {\text{H}(\R^2)}$.

			The corresponding Lie algebra defined by generators $v_1,\rho_1,v_2,\rho_2$ and relations
			$$
			[v_1,v_2]=-\rho_2, \ \ \ [\rho_1,v_2]=v_2, \ \ \ [\rho_1,\rho_2]=\rho_2.
			$$
			The complex structure is defined by
			$$
			I(v_1)=v_2, \ \ \ I(\rho_1) =\rho_2.
			$$
			For any $c,t\in\R^{>0}$ we have the following l.c.K. form  
			$$
			\omega_{c,t} =\frac{1}{c}v^1\wedge v^2+t\rho^1\wedge \rho^2
			$$
			In particular the form $\omega_{1,1}= v^1\wedge v^2+\rho^1\wedge \rho^2$ is K\"ahler.
			\end{example}
		
	


	\begin{example}
		There is an identifications $\text{SU}(2)\simeq S^3$ and a homogeneous statistical structure of curvature 1 on $S^3$. 
		The corresponding l.c.K Lie group $\left(\text{SU}(2)\times\R^{>0}\right)\ltimes \mathbb{C}^2$ is equal to the group of homothetic complex motions $\text{H}\left(\CC^2\right)$. 
		
		The algebra $\mathfrak{su}(2)$ is defined by generators $u,v,w$ and relations 
		$$
		[u,v]=2w, \ \ \ [v,w]=2u, \ \ \ [w,u]=2v.
		$$ 
		The statistical structure $\left(g,D\right)$ of constant curvature 1 on $\mathfrak{su}(2)$ is defined by 
		$$
		D_u u =D_v v = D_w w=0 , \ \ \ D_u v= w, \ \ \ D_v w = u, \ \ \ D_w u =v
		$$
		and 
		$$
		g=\left(u^*\right)^2+\left(v^*\right)^2+\left(w^*\right)^2.
		$$
		The corresponding l.c.K Lie algebra is defined by generators $u_1,v_1,w_1,\rho_1,u_2,v_2,w_2,\rho$ and relations
			$$
		[u_1,v_1]=2w_1, \ \ \ [v_1,w_1]=2u_1, \  \ \ [w_1,u_1]=2v_1, \  \ \ [u_1,v_2]=w_2, \ \ \  [v_1,w_2]=u_2, \ \ \ [w_1,u_2]=v_2,
		$$
		$$
		 [u_1,u_2]=[v_1,v_2]=[w_1,w_2]=-\rho_2, \ \ \ [\rho_1,v_2]=v_2, \ \ \ [\rho_1,u_2]=u_2, \ \ \ [\rho_1,v_2]=v_2, \ \ \  [\rho_1,\rho_2]=\rho_2.
		$$ 
		For any $c,t\in\R^{>0}$ we have the following l.c.K. form  
		$$
		\omega_{c,t} =\frac{1}{c}\left(u^1\wedge u^2+v^1\wedge v^2+w^1\wedge w^2\right)+t\rho^1\wedge \rho^2.
		$$
		In particular, the form $\omega_{1,1}=u^1\wedge u^2+v^1\wedge v^2+w^1\wedge w^2+\rho^1\wedge \rho^2$ is K\"ahler.
	\end{example}


\end{document}